\documentclass[UTF8]{article}

\usepackage{amsthm,amsmath,amssymb}
\usepackage{graphicx}
\usepackage{color}
\usepackage{ulem}

\newtheorem{theorem}{Theorem}

\newtheorem{corollary}{Corollary}

\newtheorem{definition}{Definition}
\newtheorem{proposition}{Proposition}
\newtheorem{remark}{Remark}

\newcommand{\ls}[1]
    {\dimen0=\fontdimen6\the\font\lineskip=#1\dimen0
     \advance\lineskip.5\fontdimen5\the\font
     \advance\lineskip-\dimen0
     \lineskiplimit=0.9\lineskip
     \baselineskip=\lineskip
     \advance\baselineskip\dimen0
     \normallineskip\lineskip\normallineskiplimit\lineskiplimit
     \normalbaselineskip\baselineskip
     \ignorespaces}

\UseRawInputEncoding
\begin{document}

\bibliographystyle{abbrv}

\title{On $L^0$-convex compactness in random locally convex modules}
\author{Mingzhi Wu$^1$\quad Xiaolin Zeng$^{2}$ \quad Shien Zhao$^{3}$\\
1. School of Mathematics and Physics, China University of Geosciences,\\ Wuhan {\rm 430074}, China \\
Email: wumz@cug.edu.cn\\
2. College of Mathematics and Statistics, Chongqing Technology and Business University,\\
Chongqing {\rm 400067}, China\\
Email: xlinzeng@163.com\\
3. Elementary Educational College, Capital Normal University, Beijing {\rm 100048}, China\\
Email: zsefh@cnu.edu.cn\\
}

\date{}
 \maketitle

\thispagestyle{plain}
\setcounter{page}{1}

\begin{abstract}
For the study of some typical problems in finance and economics, \v{Z}itkovi\'{c}
introduced convex compactness and gave many remarkable applications. Recently, motivated by random convex optimization and random variational inequalities, Guo, et al introduced $L^0$-convex compactness, developed the related theory of $L^0$-convex compactness in random normed modules and further applied it to backward stochastic equations.
In this paper, we extensively study $L^0$-convexly compact sets in random locally convex modules so that a series of fundamental results are obtained. First, we show that every $L^0$-convexly compact set is complete (hence is also closed and has the countable concatenation property). Then, we prove that any $L^0$-convexly compact set is linearly homeomorphic to a weakly compact subset of some locally convex space, and simultaneously establish the equivalence between $L^0$-convex compactness and convex compactness
for a closed $L^0$-convex set. Finally, we establish Tychonoff type, James type and Banach-Alaoglu type theorems for $L^0$-convex compactness, respectively.
\end{abstract}



\ls{1.5}

\section{Introduction}

Compactness plays an important role in various topics, such as optimization, variational inequalities, fixed point theory and general equilibrium theory. In an infinite-dimensional topological vector space, compactness is usually too restrictive to be used, and weak compactness often works if the space is locally convex. However, as pointed out in \cite{Zit}, quite a few optimization and equilibrium problems in mathematical economics and finance should be discussed in a non-locally convex space. To provide a proper substitute for the notion of weak compactness in a non-locally convex space, \v{Z}itkovi\'{c} \cite{Zit} introduced the notion of convex compactness as follows:

\begin{definition}(see \cite{Zit}.)
A convex subset $C$ of a topological vector space $X$ is said to be convexly compact if any family of closed and convex subsets of $C$ has a nonempty intersection whenever the family has the finite intersection property.
\end{definition}

Based on the notion of convex compactness, \v{Z}itkovi\'{c} successfully generalized some classical theorems in optimization, fixed point theory and general equilibrium theory to a more general setting, see \cite[Section 4]{Zit} for details.

Recently, random convex analysis has been developed to provide an analytic foundation for the module approach to conditional risk measures, see \cite{Guo17,Guo14,RCA1,RCA2}. According to \cite{RCA1}, random convex analysis is convex analysis on random normed modules and random locally convex modules (see Definition \ref{RLCDef} in Section 2), which are a random generalization of normed spaces and locally convex spaces, respectively. Owing to Guo's
 contribution in the last 30 years, a deep and systematic theory of random normed modules and random locally convex modules has been developed, see \cite{Guo-JFA} and the references therein for an account of main results. Aiming at the applications of random convex analysis to conditional risk measures, Guo, et al \cite{GZWW} began to
 establish some general principles for random convex optimization and random variational inequalities. The starting point of the work \cite{GZWW}
is the notion of $L^0$-convex compactness for an $L^0$-convex subset of a random locally convex module.

Throughout this paper, $(\Omega, \mathcal{F},P)$ always denotes a given probability space, $L^0(\mathcal{F}, K)$ the algebra of equivalence classes of $K$-valued random variables defined on $(\Omega,{\mathcal F},P)$, where $K$ is the scalar field ${\mathbb R}$ of real numbers or ${\mathbb C}$ of complex numbers. A subset $G$ of a random locally convex module is said to be $L^0$-convex \cite{GXC}, if $\xi x + (1-\xi)y \in G$ for any $x$ and $y$ in $G$ and any $\xi$ in $L^0({\mathcal F},{\mathbb R})$ such that $0 \leq \xi \leq 1$. We always assume that a random locally convex module is endowed with its $(\varepsilon,\lambda)$-topology (see Definition \ref{top} in Section 2).

\begin{definition}(see \cite{GZWW}.)\label{de1}
Let $E$ be a random locally convex module and $G$ a closed $L^0$-convex subset of $E$. $G$ is $L^0$-convexly compact (or, is said to have $L^0$-convex compactness), if any family of
closed $L^0$-convex subsets of $G$ has a nonempty intersection whenever this family has the finite intersection property.
\end{definition}

Using $L^0$-convex compactness, Guo, et al \cite{GZWW} successfully gave some basic results on random convex optimization and random variational inequalities in the context of a random reflexive complete random normed module, and Guo, et al \cite{GZWG} further established
the Browder-Kirk fixed point theorem in a complete random normed module and applied it to backward stochastic equations.

Since the work in \cite{GZWG,GZWW} was mainly concerned with $L^0$-convex compactness in random normed modules, in this paper we carry out an extensive study on $L^0$-convex compactness in random locally convex modules, which are more general than random normed modules. Noting the work about a new formulation of convex compactness in \cite{CDN}, we adopt a slightly different formulation of $L^0$-convex compactness as follows.

\begin{definition}\label{de2}
Let $E$ be a random locally convex module and $G$ an $L^0$-convex subset of $E$. $G$ is $L^0$-convexly compact (or, is said to have $L^0$-convex compactness), if any family of
$L^0$-convex and relatively closed subsets of $G$ has a nonempty intersection whenever this family has the finite intersection property.
\end{definition}

According to the definition, $L^0$-convex compactness is well defined only for an $L^0$-convex subset, so in the sequel when we speak of $L^0$-convex compactness of a set, we always mean that the set is $L^0$-convex.

Compared with Definition \ref{de1}, closedness of $G$ is not a premise condition in Definition \ref{de2}. Since a relatively closed subset of a closed set is closed, we see that for a closed $L^0$-convex subset of a random locally convex module, it is $L^0$-convexly compact under Definition \ref{de2} if and only if it is $L^0$-convexly compact under Definition \ref{de1}.
Surprisingly, we will show that in a random locally convex module a subset being $L^0$-convexly compact under Definition \ref{de2} is automatically closed (see Proposition \ref{ccclosed})! Thus in the context of random locally convex modules, Definition \ref{de1} and Definition \ref{de2} are in fact equivalent. Nevertheless, we insist on using Definition \ref{de2} in this paper. Following is the most important reason.

As is well-known, compactness is preserved by any continuous mapping (namely, the image of a compact space under a continuous mapping is compact), it is reasonable to require that $L^0$-convex compactness be preserved by any continuous mapping that preserves $L^0$-convexity. However, even a continuous module homomorphism does not necessarily preserve closedness, consequently, it is not easy to check whether $L^0$-convex compactness under Definition \ref{de1} is preserved or not by even a quite simple mapping (see, e.g., the proof of \cite[Lemma 2.19]{GZWW}). Comparatively, it is quite apparent that $L^0$-convex compactness under Definition \ref{de2} is preserved by some typical mappings. Precisely, we have the following self-evident proposition:

\begin{proposition}\label{lccpreserving}
Let $E_1$ and $E_2$ be two random locally convex modules over $K$ with base $(\Omega, \mathcal{F},P)$, and $G$ an $L^0$-convexly
compact subset of $E_1$. Then for any continuous mapping $f: E_1\to E_2$ such that $f(G)$ is $L^0$-convex and $f^{-1}(H)$ is $L^0$-convex for any $L^0$-convex subset $H$ of $E_2$, $f(G)$ is an $L^0$-convexly compact subset of $E_2$. Specially, for any continuous module homomorphism $f: E_1\to E_2$, $f(G)$ is an $L^0$-convexly compact subset of $E_2$.
\end{proposition}

 By Proposition \ref{lccpreserving}, we can choose a proper continuous module homomorphism to map an $L^0$-convexly compact subset of a random locally convex module to an $L^0$-convexly compact subset of some random normed module, so that we can convert $L^0$-convex compactness problems in random locally convex modules to those in random normed modules, which are easier and have been studied in depth \cite{GZWG,GZWW}.

The remainder of this paper is organized as follows: Section 2 recapitulates some basic notions and provides some simple facts. Section 3 proves that every $L^0$-convexly compact set is almost surely bounded, closed and in particular complete. Section 4 shows that any $L^0$-convexly compact set is linearly homeomorphic to a weakly compact subset of some locally convex space, and simultaneously establishes the equivalence between $L^0$-convex compactness and convex compactness for a closed $L^0$-convex sets. Section 5 presents Tychonoff type, James type and Banach-Alaoglu type theorems for $L^0$-convex compactness, respectively.

 \section{Preliminaries}

 In the sequel, we write $L^0(\mathcal{F})$ for $L^0({\mathcal F},{\mathbb R})$. Besides, $\bar{L}^0(\mathcal{F})$ denotes the set of equivalence classes of extended real-valued random variables on $(\Omega, \mathcal{F},P)$.
$\bar{L}^0(\mathcal{F})$ is usually partially ordered via $\xi\leq \eta$ if $\xi^0(\omega)\leq \eta^0(\omega)$ for almost all $\omega$ in $\Omega$, where $\xi^0$ and $\eta^0$ are arbitrarily chosen representatives of $\xi$ and $\eta$, respectively.
It is well-known from \cite{{DS58}} that $(\bar L^0({\mathcal F}), \leq)$ is a complete lattice. For any nonempty subset $H$ of $\bar L^0({\mathcal F})$, $\vee H$ and $\wedge H$ denote the supremum and infimum of $H$, respectively. As a sublattice of $\bar L^0({\mathcal F})$, $L^0({\mathcal F})$ is conditionally complete, namely any nonempty subset with an upper (resp., a lower) bound has a supremum (resp., an infimum).

Denote $L^0_+({\mathcal F})=\{ \xi \in L^0({\mathcal F})~|~ \xi \geq 0 \}$ and $L^0_{++}({\mathcal F})=\{\eta\in L^0(\mathcal F): P\{\eta>0\}=1\}$.

\begin{definition} (see \cite{Guo-JFA}.)\label{RLCDef}
 An ordered pair $(E,\mathcal{P})$ is called a random locally convex module (briefly, an $RLC$ module)
over $K$ with base $(\Omega,\mathcal{F},P)$ if $E$ is a left module over $L^0({\mathcal F},K)$ (briefly, an $L^0({\mathcal F},K)$-module) and $\mathcal{P}$ a family of mappings from $E$ to $L^0_+({\mathcal F})$ such that:\\
(RLCM-1). $\vee\{\|x\|: \|\cdot\|\in \mathcal{P}\}=0$ iff $x=\theta$ (the null element of $E$);\\
in addition, each $\|\cdot\|\in \mathcal{P}$ satisfies the following two conditions:\\
(RLCM-2). $\| \xi x \| = |\xi| \|x\|,\forall \xi \in L^0(\mathcal{F}, K), x \in E$;\\
(RLCM-3). $\|x+y\| \leq \|x\| + \|y\|, \forall x,y \in E$.\\
Furthermore, a mapping $\| \cdot \|: E\to L^0_+({\mathcal F})$ satisfying (RLCM-2) and (RLCM-3) is called an $L^0$-seminorm, in
addition, if $\|x\|=0$ also implies $x =\theta$, then it is called an $L^0$-norm, in which case $(E,\| \cdot \|)$ is
called a random normed module (briefly, an $RN$ module) over $K$ with base $(\Omega,\mathcal{F},P)$, and is a
special case of a random locally convex module when ${\mathcal P}$ consists of a single $L^0$-norm $\| \cdot \|$.
\end{definition}

When $(\Omega, \mathcal{F}, P)$ is trivial, namely $\mathcal{F}=\{\Omega, \emptyset\}$, an $RN$ or $RLC$ module over $K$ with base $(\Omega, \mathcal{F}, P)$ just reduces to an ordinary normed or locally convex space over $K$, and thus an $RN$ or $RLC$ module is a random generalization of an ordinary normed or locally convex space. The simplest $RN$ module is $(L^0(\mathcal{F},K),|\cdot|)$, where $|\cdot|$ is the absolute value mapping. It is well known that $L^0(\mathcal{F}, K)$ is a metrizable linear topological space (in fact, also a topological algebra) in the topology of convergence in probability. The $(\varepsilon,\lambda)$-topology below is exactly a natural generalization of the topology of convergence in probability.

\begin{definition}(see \cite{Guo-JFA}.)\label{top}
Let $(E, \mathcal{P})$ be an $RLC$ module over $K$ with base $(\Omega, \mathcal{F}, P)$. Denote by $\mathcal{P}_{f}$ the family of nonempty finite subsets of $\mathcal{P}$, for each $Q\in \mathcal{P}_{f}$, the $L^0$-seminorm $\|\cdot\|_Q$ is defined by $\|x\|_Q=\vee \{\|x\|: \|\cdot\| \in Q\}$ for any $x\in E$. For any given positive numbers $\varepsilon$ and $\lambda$ with $0<\lambda <1$, and for any $Q\in \mathcal{P}_{f}$, let $\mathcal{U}_{\theta}(Q,\varepsilon,\lambda)=\{x\in E: P\{\omega\in \Omega~|~ \|x\|_Q(\omega)< \varepsilon\}>1-\lambda\}$, then $\{\mathcal{U}_{\theta}(Q,\varepsilon,\lambda): Q\in \mathcal{P}_{f}, \varepsilon>0, 0<\lambda <1\}$ forms a local base of some Hausdorff linear topology for $E$, called the $(\varepsilon,\lambda)$-topology induced by $\mathcal{P}$, denoted by $\mathcal{T}_{\varepsilon,\lambda}$. Furthermore, $(E, \mathcal{T}_{\varepsilon,\lambda})$ is a Hausdorff topological module over the topological algebra $L^0(\mathcal{F}, K)$.
\end{definition}

It is easy to see that in an $RLC$ module $(E, \mathcal{P})$ over $K$ with base $(\Omega, \mathcal{F}, P)$, a net $\{x_{\alpha}, \alpha\in \Gamma\}$ converges in $\mathcal{T}_{\varepsilon,\lambda}$ to $x$ iff $\{\|x_{\alpha}-x\|, \alpha\in \Gamma\}$ converges in probability to 0 for each $\|\cdot\|\in \mathcal{P}$. In particular, in an $RN$ module $(E, \|\cdot\|)$, $\mathcal{T}_{\varepsilon,\lambda}$ is metrizable, and a sequence $\{x_n, n\in \mathbb{N}\}$ converges in $\mathcal{T}_{\varepsilon,\lambda}$ to $x$ iff $\{\|x_n-x\|, n\in \mathbb{N}\}$ converges in probability to 0.

Any $RLC$ module can be embedded densely into a complete $RLC$ module, for later use, we restate \cite[Lemma 3.19]{Guo-JFA} as follows.
\begin{proposition}\label{completion}
Let $(E,\mathcal{P})$ be an $RLC$ module over $K$ with base $(\Omega,\mathcal{F},P)$. Then there exists an $L^0(\mathcal{F},K)$-module $\tilde E\supset E$ such that each $\|\cdot\|\in {\mathcal P}$ can be extended to an $L^0$-seminorm on $\tilde E$ in a unique way (still denoted by $\|\cdot\|$), so that $(\tilde E,\mathcal{P})$ is a complete $RLC$ module and $E$ is dense in $\tilde E$. Besides, such $(\tilde E,\mathcal{P})$ is unique in the sense of $\mathcal{P}$-isometric isomorphism, called the completion of $(E,\mathcal{P})$.
\end{proposition}

In this paper, we need to consider the topological product space of a family of $RLC$ modules. We give an obvious fact as follows.
\begin{proposition}
Let $\{(E_i,\mathcal{P}_i), i\in I\}$ be a family of $RLC$ modules over $K$ with base $(\Omega,\mathcal{F},P)$. Let $E=\prod_{i\in I}E_i$, then in an apparent way $E$ becomes an $L^0(\mathcal{F},K)$-module. Fixing $i\in I$, $\|\cdot\|\circ \pi_i$ is an $L^0$-seminorm on $E$ for any $\|\cdot\|\in {\mathcal P}_i$, where $\pi_i:E\to E_i$ is the canonical projection. Denote ${\mathcal P}=\{\|\cdot\|\circ \pi_i: i\in I, \|\cdot\|\in \mathcal{P}_i\}$, then $(E, {\mathcal P})$ is an $RLC$ module and its $(\varepsilon,\lambda)$-topology is just the product topology for $\prod_{i\in I}(E_i,\mathcal{P}_i)$. Thus the topological product space of a family of $RLC$ modules is also an $RLC$ module.
\end{proposition}

\section{Almost Surely Boundedness, Closedness and Completeness}

As is well-known, a weakly compact subset of a locally convex space is bounded, closed and complete. In this section, we will show that every $L^0$-convexly compact subset of an $RLC$ module is almost surely bounded, closed and in particular complete.

We start with closedness.

Proposition \ref{ccclosed} below is similar to \cite[Lemma 1]{CDN}, which says that a convexly compact subset of a locally convex space must be closed. The proof of \cite[Lemma 1]{CDN} is based on the definition that a locally convex space has a local base whose members are convex, however for an $RLC$ module endowed with the $(\varepsilon,\lambda)$-topology, there does not exist any local base whose members are $L^0$-convex,  as a result we can not give a proof of Proposition \ref{ccclosed} by directly following the proof of \cite[Lemma 1]{CDN}. This observation forces us to use the locally $L^0$-convex topology for an $RLC$ module. The locally $L^0$-convex topology was introduced in \cite{FKV09}, since this paper involves this topology only once, to save space, we do not restate it.

As pointed out in \cite{Guo-JFA}, the bridge between the $(\varepsilon,\lambda)$-topology and the locally $L^0$-convex topology is the notion of the countable concatenation property.
For the reader's convenience, let us recall \cite{Guo-JFA}: let $G$ be a subset of an $L^0(\mathcal{F}, K)$-module $E$, if for any sequence $\{ g_n, n \in \mathbb{N} \}$ in $G$ and any countable partition $\{ A_n: n \in \mathbb{N}\}$ of $\Omega$ to $\mathcal{F}$, there is $g \in G$ such that $\tilde{I}_{A_n} g = \tilde{I}_{A_n} g_n$ for every $n \in \mathbb{N}$, then we say that $G$ has the countable concatenation property, where $\tilde{I}_A$ denotes the equivalence class of the characteristic function $I_A$ for any $A \in \mathcal{F}$.

We first show that every $L^0$-convexly compact subset of an $RLC$ module must have the countable concatenation property.

In the sequel, for a subset $G$ of an $RLC$-module, $conv_{L^0}(G)$ denotes the $L^0$-convex hull of $G$ (namely, the smallest $L^0$-convex set containing $G$), and $\overline {G}$ the closure of $G$.

\begin{proposition}\label{ccp}
 Let $(E,{\mathcal P})$ be an $RLC$ module over $K$ with base $(\Omega,{\mathcal F},P)$
 and $G$ an $L^0$-convexly compact subset of $E$, then $G$ must have the countable concatenation property.
\end{proposition}

\begin{proof}
We can, without loss of generality, assume that the null element $\theta$ of $E$ belongs to $G$, otherwise, we make a translation.
  Let $\{x_n, n\in \mathbb{N}\}$ be an arbitrary sequence in $G$ and $\{A_n: n\in \mathbb{N}\}$ an arbitrary countable partition of $\Omega$ to ${\mathcal F}$. It is enough to show that there exists an $x\in G$ such that ${\tilde I}_{A_n}x={\tilde I}_{A_n}x_n$ for every $n\in \mathbb{N}$.
For each $n\in \mathbb{N}$, let $y_n={\tilde I}_{A_1}x_1+{\tilde I}_{A_2}x_2+\cdots+{\tilde I}_{A_n}x_n$, by the $L^0$-convexity of $G$ and the assumption $\theta\in G$ we see that $y_n\in G$. Further, let $G_n=\overline {conv_{L^0}\{y_k,k\geq n\}}\cap G$ for each $n\in \mathbb{N}$, then $\{G_n, n\in \mathbb{N}\}$ is a sequence of $L^0$-convex and relatively closed subset of $G$ with the finite intersection property. By the construction, each $z$ in ${conv_{L^0}\{y_k,k\geq n\}}$ must satisfy ${\tilde I}_{A_k}z={\tilde I}_{A_k}x_k$ for every $k\in \{1,2,\dots,n\}$, which in turn implies that each $z$ in $\overline {conv_{L^0}\{y_k,k\geq n\}}$ possesses the same property. By the $L^0$-convex compactness of $G$, there exists some $x$ in $E$ such that $x$ belongs to every $G_n$. We conclude that $x$ belongs to $G$ and ${\tilde I}_{A_k}x={\tilde I}_{A_k}x_k$ for every $k\in \mathbb{N}$. \quad$\square$
\end{proof}

Now we can show that every $L^0$-convexly compact subset of an $RLC$ module must be closed.

\begin{proposition}\label{ccclosed}
 Let $(E,{\mathcal P})$ be an $RLC$ module over $K$ with base $(\Omega,{\mathcal F},P)$
 and $G$ an $L^0$-convexly compact subset of $E$, then $G$ must be closed.
\end{proposition}
\begin{proof}
We prove $G$'s closedness by contradiction. Suppose that $G$ is not closed and let $x_0$ be an element in ${\overline G}\setminus G$. Since $G$ is $L^0$-convexly compact, $G$ has the countable concatenation property by Proposition \ref{ccp}, then it follows from \cite[Theorem 3.12]{Guo-JFA} that ${\overline G}={\overline G}_c$, where ${\overline G}_c$ is the closure of $G$ under the locally $L^0$-convex topology. Let $\mathcal{P}_{f}$ denote the family of nonempty finite subsets of ${\mathcal P}$. For any $\mathcal Q\in \mathcal{P}_{f}$ and any $\xi\in L^0_{++}({\mathcal F})$, let $U_\theta({\mathcal Q, \xi})=\{x\in E:\|x\|_{\mathcal Q}\leq \xi\}$, where the $L^0$-seminorm $\|\cdot\|_{\mathcal Q}$ is defined as in Definition \ref{top}, then by \cite[Proposition 2.7]{Guo-JFA}, $\{U_\theta({\mathcal Q, \xi}): \mathcal Q\in \mathcal{P}_{f},~\xi\in L^0_{++}({\mathcal F})\}$ is a local base of the locally $L^0$-convex topology for $E$. Since $x_0\in {\overline G}_c$, each $F_{\mathcal Q, \xi}:=(x_0+U_\theta({\mathcal Q, \xi}))\cap G$ is nonempty. Given $Q_i\in \mathcal{P}_{f}$ and $\xi_i\in L^0_{++}({\mathcal F})$ for each $i\in\{1,2,\dots,d\}$, let $\mathcal Q=\mathcal Q_1\cup\cdots\cup\mathcal Q_d$ and $\xi=\wedge\{\xi_1,\dots,\xi_d\}$, then $\mathcal Q\in \mathcal{P}_{f}$ and $\xi\in L^0_{++}({\mathcal F})$ so that $U_\theta({\mathcal Q, \xi})\subset U_\theta({\mathcal Q_i, \xi_i}),i=1,\dots,d$, which in turn implies that $F_{\mathcal Q, \xi}\subset F_{\mathcal Q_1, \xi_1}\cap\cdots\cap F_{\mathcal Q_d, \xi_d}$. Thus the family $\{F_{\mathcal Q, \xi}: \mathcal Q\in \mathcal{P}_{f},~\xi\in L^0_{++}({\mathcal F})\}$ has the finite intersection property. Now each $U_\theta({\mathcal Q, \xi})$ is $L^0$-convex and closed under the $(\varepsilon,\lambda)$-topology, it follows that each $F_{\mathcal Q, \xi}$ is an $L^0$-convex and relatively closed subset of $G$, then by the $L^0$-convex compactness of $G$, we obtain that $\bigcap\{F_{\mathcal Q, \xi}: \mathcal Q\in \mathcal{P}_{f}, \xi\in L^0_{++}({\mathcal F})\}=\{x_0\}\cap G$ is nonempty. This contradicts to the assumption $x_0\notin G$. \quad$\square$
\end{proof}

We then turn to almost surely boundedness.

In line with \cite{GZWW}, a nonempty subset $G$ of an $RLC$ module $(E,{\mathcal P})$ is said to be almost surely bounded if for each $\|\cdot\|\in {\mathcal P}$, there is some $\xi\in L^0_+({\mathcal F})$ such that $\|x\|\leq \xi, \forall x\in G$. It has been shown in \cite[Lemma 2.19]{GZWW} that a closed $L^0$-convexly compact subset of an $RLC$ module is almost surely bounded. In view of Proposition \ref{ccclosed}, in an $RLC$ module, an $L^0$-convexly compact set under Definition \ref{de2} is just a closed $L^0$-convexly compact set under Definition \ref{de1}, we thus obtain the following:

\begin{proposition}\label{asb}
 Let $(E,{\mathcal P})$ be an $RLC$ module over $K$ with base $(\Omega,{\mathcal F},P)$
 and $G$ an $L^0$-convexly compact subset of $E$, then $G$ must be almost surely bounded.
\end{proposition}

In the end of this section, we show that an $L^0$-convexly compact set is complete. We point out that completeness is the most desired property in this paper.

\begin{theorem}\label{complete}
  Let $(E,{\mathcal P})$ be an $RLC$ module over $K$ with base $(\Omega,{\mathcal F},P)$
 and $G$ an $L^0$-convexly compact subset of $E$, then $G$ is complete.
\end{theorem}
\begin{proof}
Let $(\tilde E,{\mathcal P})$ be the completion of $(E,{\mathcal P})$ (see Proposition \ref{completion}) and $i: E\to \tilde E$ the inclusion mapping. Since $i$ is a continuous module homomorphism, it follows from Proposition \ref{lccpreserving} that $G=i(G)$ is $L^0$-convexly compact in $(\tilde E,{\mathcal P})$. Then $G$ is closed in $(\tilde E,{\mathcal P})$ by Proposition \ref{ccclosed}. As a closed subset of the complete $RLC$ module $(\tilde E,{\mathcal P})$, $G$ is thus complete. \quad$\square$
\end{proof}

\section{$L^0$-convex compactness VS. Convex compactness}

In the proof of \cite[Theorem 2.21]{GZWW}, it has been shown that an $L^0$-convexly compact subset of a complete $RN$ module is linearly homeomorphic to a convexly compact subset of some Banach space, thus for a closed $L^0$-convex subset of a complete $RN$ module, $L^0$-convex compactness is equivalent to convex compactness. In this section, we will extend this equivalence to the context of an $RLC$ module.

We first give an embedding theorem as follows.

\begin{proposition}\label{embed}
Let $(E,{\mathcal P})$ be an $RLC$ module over $K$ with base $(\Omega,{\mathcal F},P)$, then there exists a family of $RN$ modules $\{(E_p,\|\cdot\|_p), p\in {\mathcal P}\}$ over $K$ with base $(\Omega,{\mathcal F},P)$ together with a topological module homeomorphism $h$ that carries $E$ onto an $L^0({\mathcal F},K)$-submodule of the product space $\prod_{p\in {\mathcal P}}(E_p, \|\cdot\|_p)$. Besides, we can require that each $(E_p,\|\cdot\|_p)$ be complete.
\end{proposition}

\begin{proof}
 For each $p\in {\mathcal P}$, denote $N_p=\{x\in E:p(x)=0\}$. Since $p$ is a continuous $L^0$-seminorm on $E$, $N_p$ is a closed $L^0({\mathcal F},K)$-submodule of $(E,{\mathcal T}_{\varepsilon,\lambda})$. Let $E_p=E/N_p$ be the quotient module and $\pi_p: E\to E_p$ the canonical quotient mapping. Noting that $p(x)=p(y)$ whenever $x$ and $y$ are any two elements of $E$ such that $\pi_p(x)=\pi_p(y)$, we can define $\|\cdot\|_p: E_p\to L^0_+({\mathcal F})$ by $\|\pi_p(x)\|_p=p(x)$ for every $x\in E$. It is obvious that $\|\cdot\|_p$ is an $L^0$-norm on $E_p$ and $\pi_p: (E,{\mathcal P})\to (E_p,\|\cdot\|_p)$ is a continuous module homomorphism.

 Define $h: (E,{\mathcal P})\to \prod_{p\in {\mathcal P}}(E_p, \|\cdot\|_p)$ by $h(x)=(\pi_p(x))_{p\in {\mathcal P}}$ for every $x\in E$. Since each $\pi_p$ is a module homomorphism, it is clear that $h$ is also a module homomorphism, so that $h(E)$ is an $L^0({\mathcal F},K)$-module. It follows from (RLCM-1) that $h$ is injective. According to the definition of the $(\varepsilon,\lambda)$-topology, it is straightforward to check that $h: (E,{\mathcal T}_{\varepsilon,\lambda})\to (h(E),{\mathcal T}_{\varepsilon,\lambda})$ is a topological module homeomorphism.

If some $RN$ module $(E_p, \|\cdot\|_p)$ is not complete, then we can substitute it with its completion (see Proposition \ref{completion}), thus we can require that each $(E_p,\|\cdot\|_p)$ be complete. \quad$\square$
\end{proof}

Proposition \ref{embed} provides a way to convert a problem in an $RLC$ module to that in an $RN$ module. However, we need further convert a problem in an $RN$ module to that in an ordinary normed space.
In this aspect, a powerful tool is the abstract $L^2$ space generated from an $RN$ module: let $(E,\|\cdot\|)$ be an $RN$ module over $K$ with base $(\Omega,{\mathcal F},P)$, denote $L^2(E)=\{x\in E~|~\int_\Omega \|x\|^2 dP<+\infty\}$, then $(L^2(E),\|\cdot\|_2)$ is a normed space over $K$, where the norm $\|\cdot\|_2$ is given by $\|x\|_2=(\int_\Omega \|x\|^2 dP)^{\frac{1}{2}}$ for any $x\in L^2(E)$, furthermore, if $(E,\|\cdot\|)$ is complete, then $(L^2(E),\|\cdot\|_2)$ is a Banach space. The abstract $L^2$ space generated from an $RN$ module has been frequently used in the study of $RN$ modules and $RLC$ modules, see \cite{GL,GXC,GZ,GZWW} for examples.

In the sequel, given two linear spaces $X_1$ and $X_2$ (over $K$), a mapping $T$ from a convex subset $G_1$ of $X_1$ to a convex subset $G_2$ of $X_2$ is said to be linear if $T$ is a restriction of a linear mapping from $X_1$ to $X_2$.

\begin{theorem}\label{wchomo}
Let $(E,{\mathcal P})$ be an $RLC$ module over $K$ with base $(\Omega,{\mathcal F},P)$ and $G$ an $L^0$-convex subset of $E$. Then the following two statements are equivalent:\\
(1). $G$ is $L^0$-convexly compact.\\
(2). There exists a complete locally convex space $(X,\mathcal T)$ over $K$ together with a linear mapping $f:G\to X$ such that $(G,{\mathcal T}_{\varepsilon,\lambda})$ is homeomorphic to $(f(G), \mathcal T)$ under $f$ and $f(G)$ is a convex and weakly compact subset of $X$.
\end{theorem}

\begin{proof}
$(2)\Longrightarrow (1)$. Suppose that $\{G_a, a\in A\}$ is a family of $L^0$-convex and relatively closed subset of $G$ with the finite intersection property, since $f: (G,{\mathcal T}_{\varepsilon,\lambda})\to (f(G), \mathcal T)$ is a linear homeomorphism, $\{f(G_a), a\in A\}$ is a family of convex and relatively closed subsets of $(f(G), \mathcal T)$ with the finite intersection property. Since a convex relatively closed subset is also relatively closed in the weak topology, it follows from weak compactness of $f(G)$ that $\{f(G_a), a\in A\}$ has nonempty intersection. Thus $\{G_a, a\in A\}$ has nonempty intersection because $f$ is injective. Therefore, $G$ is $L^0$-convexly compact.

$(1)\Longrightarrow (2)$. Assume that $G$ is $L^0$-convexly compact.
For each $p\in {\mathcal P}$, let $(E_p, \|\cdot\|_p)$ and $\pi_p$ be the same as in the proof of Proposition \ref{embed} and assume that $(E_p, \|\cdot\|_p)$ is complete. Since $\pi_p: (E,{\mathcal P})\to (E_p,\|\cdot\|_p)$ is a continuous module homomorphism, it follows from Proposition \ref{lccpreserving} that $\pi_p(G)$ is $L^0$-convexly compact in $(E_p, \|\cdot\|_p)$, thus by Proposition \ref{asb} there exists $\xi_p\in L^0_{++}({\mathcal F})$ such that $\|\pi_p(x)\|_p\leq \xi_p$ for any $x\in G$. Let $G_p=\frac{\pi_p(G)}{\xi_p}$, then $\|x\|_p\leq 1$ for every $x\in G_p$, thus $G_p\subset L^2(E_p)$. Moreover, as shown in the proof of \cite[Theorem 2.21]{GZWW}, $G_p$ is a weakly compact subset of the Banach space $(L^2(E_p),\|\cdot\|_2)$. If $w$ denotes the weak topology, then by Tychonoff's theorem, $\prod_{p\in {\mathcal P}} G_p$ is a compact subset of the topological product space $\prod_{p\in {\mathcal P}}(L^2(E_p),w)$.

Making a slight modification of $h$ in the proof of Proposition \ref{embed}, define $h: (E,{\mathcal P})\to \prod_{p\in {\mathcal P}}(E_p, \|\cdot\|_p)$ by $h(x)=(\frac{\pi_p(x)}{\xi_p})_{p\in {\mathcal P}}$ for every $x\in E$. Then it is easy to see that $h$ is still a topological module homeomorphism from $(E,{\mathcal T}_{\varepsilon,\lambda})$ to $(h(E),{\mathcal T}_{\varepsilon,\lambda})$.

Let $(X, \mathcal T)$ be the topological product space $\prod_{p\in {\mathcal P}}(L^2(E_p),\|\cdot\|_2)$ (where each factor is endowed with the norm topology). Now that each factor is a Banach space, as is well-known \cite[Theorem 7.2, p.57]{KN}, $(X,\mathcal T)$ is a complete locally convex space. Please note that since $X=\prod_{p\in {\mathcal P}}L^2(E_p)\subset \prod_{p\in {\mathcal P}}E_p$, $X$ also inherits the topology ${\mathcal T}_{\varepsilon,\lambda}$ from $\prod_{p\in {\mathcal P}}(E_p, \|\cdot\|_p)$.

 Obviously we have $h(G)\subset \prod_{p\in {\mathcal P}}G_p\subset X$. Define $f: G\to X$ as the restriction of $h$ to $G$, then $f$ is a linear homeomorphism from $(G,{\mathcal T}_{\varepsilon,\lambda})$ to $(f(G),{\mathcal T}_{\varepsilon,\lambda})$. Using Lebesgue's dominated convergence theorem, the two topologies--${\mathcal T}_{\varepsilon,\lambda}$ and ${\mathcal T}$--coincide on $f(G)$. Therefore $f$ is a linear homeomorphism from $(G,{\mathcal T}_{\varepsilon,\lambda})$ to $(f(G), \mathcal T)$.

 From Theorem \ref{complete}, $G$ is complete, then $f(G)$ is a complete and thus closed subset of $(X,\mathcal T)$. Noting that $f(G)$ is convex, $f(G)$ is thus weakly closed in $X$. As is known \cite[Theorem 17.13, p.160]{KN}, the weak topology for the product of locally convex spaces is the product of weak topologies for factor spaces, thus $(X,w)=\prod_{p\in {\mathcal P}}(L^2(E_p),w)$. We conclude that $f(G)$, as a closed subset of a compact set $\prod_{p\in {\mathcal P}} G_p$ in $(X,w)$, is also a compact subset of $(X,w)$. In other words, $f(G)$ is a weakly compact subset of $(X,\mathcal T)$. \quad$\square$
\end{proof}

\begin{remark}
In the proof of ``$(2)\Longrightarrow (1)$'', we, in fact, have shown that the statement (2) implies that $G$ is convexly compact, thus an $L^0$-convexly compact is convexly compact. Obviously, a closed and convexly compact $L^0$-convex set is also $L^0$-convexly compact, hence for a closed $L^0$-convex subset of an $RLC$ module, $L^0$-convex compactness is equivalent to convex compactness. This extends the equivalence pointed out in \cite[Remark 2.22]{GZWW} in the context of a complete $RN$ module. Compared to \cite[Remark 2.22]{GZWW}, the equivalence holds even without the completeness assumption on the $RLC$ module. In fact, the completeness assumption on $X$ in the statement (2) is also redundant, as we know that a weakly compact subset of a locally convex space is always complete (see \cite[Corollary 17.8, p.156]{KN}).
\end{remark}

\section{Tychonoff, James and Banach-Alaoglu type theorems for $L^0$-convex compactness}

Tychonoff theorem, James theorem and Banach-Alaoglu theorem are the most fundamental theorems concerning compactness, weak compactness and weak-star compactness, respectively. In this section, we establish their variants for $L^0$-convex compactness.

We first establish a Tychonoff type theorem for the topological product of a family of $L^0$-convexly compact sets.

\begin{theorem}\label{tychonoff}
 Suppose that $\{(E_i, {\mathcal P}_i), i\in I\}$ is a family of $RLC$ modules over $K$ with base $(\Omega,{\mathcal F},P)$, and $G_i$ is an $L^0$-convexly compact subset of $E_i$ for each $i\in I$. Then $G=\prod_{i\in I}G_i$ is $L^0$-convexly compact in $\prod_{i\in I}(E_i, {\mathcal P}_i)$.
\end{theorem}
\begin{proof} Obviously, $G$ is $L^0$-convex.
By Theorem \ref{wchomo}, for each $i\in I$, there exists a complete locally convex space $(X_i,{\mathcal T}_i)$ together with a linear mapping $f_i: G_i\to X_i$ such that $(G_i,{\mathcal T}_{\varepsilon,\lambda})$ is homeomorphic to $(f_i(G_i), \mathcal T_i)$ under $f_i$ and $f_i(G_i)$ is a weakly compact subset of $(X_i,{\mathcal T}_i)$. Let $(X,\mathcal T)$ be the topological product space $\prod_{i\in I}(X_i,{\mathcal T}_i)$, then from \cite[Theorem 7.2, p.57]{KN} we know that $(X,\mathcal T)$ is a complete locally convex space. We use $w$ to denote the weak topology for any locally convex space. Define $f: G\to X$ by $f(x)=(f_i(x_i))_{i\in I}$ for every $x=(x_i)_{i\in I}\in G$. Then it is straightforward to check that: (1) $f$ is linear and injective, and (2) $f: (G, {\mathcal T}_{\varepsilon,\lambda})\to (f(G), \mathcal T)$ is a homeomorphism. Since each factor $G_i$ is complete by Theorem \ref{complete}, then $G$ is complete, therefore $f(G)$ is a complete and thus closed subset of $(X,\mathcal T)$. Using the fact that $f(G)$ is convex, $f(G)$ is also a closed subset of $(X,w)=\prod_{i\in I}(X_i,w)$. By Tychonoff's theorem, $\prod_{i\in I}f_i(G_i)$ is a compact subset of $(X,w)$. Therefore, $f(G)$, as a closed subset of $\prod_{i\in I}f_i(G_i)$, is also a compact subset of $(X,w)$.
By Proposition \ref{wchomo}, $G$ is an $L^0$-convexly compact subset of $\prod_{i\in I}(E_i, {\mathcal P}_i)$. \quad$\square$
\end{proof}

 In order to establish James type and Banach-Alaoglu type theorems for $L^0$-convex compactness, we need the notions of random conjugate space and random weak-star topology. Let $(E,{\mathcal P})$ be an $RLC$ module over $K$ with base $(\Omega,{\mathcal F},P)$. Then the set of continuous module homomorphisms from $(E,{\mathcal P})$ to $L^0(\mathcal{F},K)$ forms an $L^0(\mathcal{F},K)$-module in a natural way, denoted by $E^\ast$ and called the random conjugate space of $E$. Specially, for an $RN$ module $(E,\|\cdot\|)$, a module homomorphism $f: E\to L^0(\mathcal{F},K)$ is an element of $E^\ast$ iff there exists a $\xi\in L^0_+({\mathcal F})$ such that $|f(x)|\leq \xi\|x\|, \forall x\in E$, at this time, denote $\|f\|^\ast=\wedge\{\xi\in L^0_+({\mathcal F})|~|f(x)|\leq \xi\|x\| \text {~for all~} x\in E\}$, then $\|\cdot\|^\ast$ is an $L^0$-norm on $E^\ast$ so that $(E^\ast,\|\cdot\|^\ast)$ becomes an $RN$ module. For an $RLC$ module $(E,{\mathcal P})$, each $x\in E$ induces an $L^0$-seminorm $|\langle\cdot,x\rangle|$ on $E^\ast$ by $|\langle f,x\rangle|=|f(x)|, \forall f\in E^\ast$, then $\sigma(E^\ast,E):=\{|\langle\cdot,x\rangle|: x\in E\}$ is a family of $L^0$-seminorms such that $(E^\ast, \sigma(E^\ast,E))$ is an $RLC$ module, whose $(\varepsilon,\lambda)$-topology is called the random weak-star topology. The random weak-star topology was first introduced in \cite{BBKS}.

We then present a James type characterization theorem for $L^0$-convex compactness in a complete $RLC$ module, which is an extension of \cite[Theorem 2.21]{GZWW}.

\begin{theorem}
Let $(E,{\mathcal P})$ be a complete $RLC$ module over $K$ with base $(\Omega,{\mathcal F},P)$ and $G$ a closed $L^0$-convex subset of $E$. Then $G$ is $L^0$-convexly compact if and only if for each $f\in E^\ast$ there exists $g_0\in G$ such that $Re(f(g_0))=\vee\{Re(f(g)):g\in G\}$,  where, for any $\xi\in L^0(\mathcal{F},K)$, $Re(\xi)$ denotes the real part of $\xi$.
\end{theorem}

\begin{proof}
Necessity. Clearly $Re f: E\to L^0({\mathcal F})$ is continuous and $L^0$-linear, namely $Re f(sx+ty)=sRef(x)+tRef(y)$ for any $x$ and $y$ in $E$ and any $s$ and $t$ in $L^0({\mathcal F})$. Then $Ref(G)$ is $L^0$-convex and $(Ref)^{-1}(H)$ is $L^0$-convex for any $L^0$-convex subset $H$ of $L^0({\mathcal F})$. It follows from Proposition \ref{lccpreserving} that $Re f(G)$ is an $L^0$-convexly compact subset of $(L^0({\mathcal F}),|\cdot|)$. Then by \cite[Theorem 2.17]{GZWW}, $Ref(G)$ is a random closed interval $[m, M]:=\{\xi\in L^0({\mathcal F}): m\leq \xi\leq M\}$, where $m=\wedge\{Ref(g):g\in G\}$ and $M=\vee\{Ref(g):g\in G\}$, in particular, there exists $g_0\in G$ such that $Re(f(g_0))=M=\vee\{Re(f(g)):g\in G\}$.

Sufficiency. Let $\{(E_p,\|\cdot\|_p), p\in {\mathcal P}\}$ and $h: (E,{\mathcal P})\to \prod_{p\in {\mathcal P}}(E_p, \|\cdot\|_p)$ be the same as in the proof of Proposition \ref{embed} and assume that each $(E_p, \|\cdot\|_p)$ is complete. For each $q\in {\mathcal P}$, we now write $\pi_q$ for the canonical projection from $\prod_{p\in {\mathcal P}}(E_p, \|\cdot\|_p)$ to $(E_q, \|\cdot\|_q)$, and $G_q$ the closure of $\pi_qh(G)$ in $(E_q, \|\cdot\|_q)$. Since $\pi_qh(G)$ is $L^0$-convex, so is its closure $G_q$. For any $f_q\in E^\ast_q$, $f_q\pi_qh$ belongs to $E^\ast$, thus $Re(f_q\pi_qh)$ attains its supremum on $G$, in other words, $Re(f_q)$ attains its supremum on $\pi_qh(G)$, which is equal to supremum on $G_q$, the closure of $\pi_qh(G)$. Then $G_q$ is $L^0$-convexly compact by \cite[Theorem 2.21]{GZWW}. Thus by Theorem \ref{tychonoff}, $\prod_{q\in {\mathcal P}}G_q$ is an $L^0$-convexly compact subset of $\prod_{q\in {\mathcal P}}(E_q, \|\cdot\|_q)$. Since $G$ is a closed subset of the complete $RLC$ module $(E,{\mathcal P})$, $G$ is complete, hence $h(G)$ is a complete and thus closed $L^0$-convex subset of $\prod_{q\in {\mathcal P}}G_q$. It follows that $h(G)$ is $L^0$-convexly compact, which implies that $G$ is $L^0$-convexly compact. \quad$\square$
\end{proof}

 We conclude this paper with a Banach-Alaoglu type theorem for $L^0$-convex compactness. We would like to remind the reader of the work \cite{BBKS}. Generally, the random closed unit ball of the random conjugate space of an $RN$ module is not compact under the random weak-star topology (see \cite[Theorem 2.3]{BBKS}), whereas we will show that it is always $L^0$-convexly compact under the random weak-star topology (see Corollary \ref{BARN} below).

\begin{theorem}\label{BA}
Let $(E,{\mathcal P})$ be an $RLC$ module over $K$ with base $(\Omega,{\mathcal F},P)$ and $E^\ast$ its random conjugate space. If $V$ is a subset of $E$ containing the random closed unit ball $B_p:=\{x\in E: p(x)\leq 1\}$ of some continuous $L^0$-seminorm $p$ on $E$, then its polar $V^o:=\{x^\ast\in E^\ast: |x^\ast(x)|\leq 1, \forall x\in V\}$ is an $L^0$-convexly compact subset of $(E^\ast, \sigma(E^\ast,E))$.
\end{theorem}
\begin{proof}
It is obvious that $V^o$ is an $L^0$-convex and closed subset of $(E^\ast, \sigma(E^\ast,E))$. From the assumption on $V$, there exists a $\xi_x\in L^0_{++}({\mathcal F})$ for each $x\in E$ such that $x\in \xi_xV$ (in fact, since $p(\frac{x}{p(x)+1})=\frac{p(x)}{p(x)+1}\leq 1$, we have $\frac{x}{p(x)+1}\in V$ so that $x=(p(x)+1)\frac{x}{p(x)+1}\in (p(x)+1)V$), as a result $|x^\ast(x)|\leq \xi_x$ for any $x^\ast\in V^o$. For each $x\in E$, let $I_x=\{\xi\in L^0({\mathcal F},K):|\xi|\leq \xi_x\}$, then from \cite[Corollary 2.23]{GZWW} $I_x$ is $L^0$-convexly compact in the complete random reflexive $RN$ module $(L^0({\mathcal F},K),|\cdot|)$. It follows from Theorem \ref{tychonoff} that $H:=\prod_{x\in E}I_x$ is $L^0$-convexly compact in $\prod_{x\in E}(L^0({\mathcal F},K),|\cdot|)$. As in the proof of classical Banach-Alaoglu theorem, $V^o$ is regarded as a subset of $H$, and the random weak-star topology inherited from $E^\ast$ and the topology ${\mathcal T}_{\varepsilon,\lambda}$ inherited from $H$ coincide on $V^o$. It remains to show that $V^o$ is closed in $H$. Note that $H$ is complete, suppose that $f\in H$ is the ${\mathcal T}_{\varepsilon,\lambda}$-limit of a net $\{x^\ast_\alpha,\alpha\in D\}$ in $V^o$, we only need to show that $f\in V^o$.
Clearly, $f$ is a module homomorphism from $E$ to $L^0({\mathcal F},K)$ and $|f(x)|\leq 1, \forall x\in V$. Fixing an $x\in E$, clearly $\frac{x}{p(x)+n^{-1}}\in B_p\subset V$ for any positive integer $n$, thus $|f(x)|\leq p(x)+n^{-1}$. Letting $n\to\infty$, we see that $|f(x)|\leq p(x)$. Hence we obtain that $|f(x)|\leq p(x), \forall x\in E$, which implies that $f\in E^\ast$. \quad$\square$
\end{proof}

Corollary \ref{BARN} below is the most important case of Theorem \ref{BA}.

\begin{corollary}\label{BARN}
Let $(E,\|\cdot\|)$ be an $RN$ module over $K$ with base $(\Omega,{\mathcal F},P)$ and $(E^\ast,\|\cdot\|^\ast)$ its random conjugate space. Then the random closed unit ball $\{x^\ast\in E^\ast: \|x^\ast\|^\ast\leq 1\}$ of $E^\ast$ is an $L^0$-convexly compact subset of $(E^\ast, \sigma(E^\ast,E))$.
\end{corollary}

\begin{proof}
Let $V$ be the random closed unit ball of $E$, namely $V=\{x\in E: \|x\|\leq 1\}$, then $\{x^\ast\in E^\ast: \|x^\ast\|^\ast\leq 1\}=V^o$, thus the conclusion follows from Theorem \ref{BA}. \quad$\square$
\end{proof}

\begin{remark}
In the statement of the classical Banach-Alaoglu theorem (see, e.g. \cite[Chapter 3, 3.15]{RuFA}, or \cite[Theorem 17.4, p.155]{KN}), $V$ is a neighborhood of $\theta$ of a topological vector space, whereas in the statement of our Theorem \ref{BA}, $V$ is not necessarily a ${\mathcal T}_{\varepsilon,\lambda}$-neighborhood of $\theta$ of an $RLC$ module. There are two reasons for this change. One reason is that when the base $(\Omega,{\mathcal F},P)$ is trivial (namely, ${\mathcal F}=\{\Omega,\emptyset\}$), the $RLC$ module $(E,{\mathcal P})$ reduces to a locally convex space, and the random unit ball of a continuous $L^0$-seminorm reduces to the unit ball of a continuous seminorm so that $V$ is a neighborhood of $\theta$, in this case Theorem \ref{BA} reduces to the classical Banach-Alaoglu theorem. The other more important reason is that a ${\mathcal T}_{\varepsilon,\lambda}$-neighborhood of $\theta$ of an $RLC$ module is too `` large '' to be widely used, for instance, it is easy to see that the random unit ball of an $RN$ module is not a ${\mathcal T}_{\varepsilon,\lambda}$-neighborhood of $\theta$ (unless the base $(\Omega,{\mathcal F}, P)$ is essentially generated by finitely many atoms). In view of Corollary \ref{BARN}, the formulation of Theorem \ref{BA} is proper.
\end{remark}

\section*{Acknowledgements}
 The first author is supported by the Natural Science Foundations of China (Grant No.11701531) and the Fundamental Research Funds for the Central Universities, China University of Geosciences (Wuhan)(Grant No. CUGL170820). The second author is supported by Natural Science Foundation of Chongqing, China (Grant No. cstc2020jcyj-msxmX0328) and the Science and Technology Research Program of Chongqing Municipal Education Commission (Grant No. KJQN202000838).


\end{document}